\documentstyle[12pt,amssymb,amscd,psfig]{amsart}
\newcommand{\np}{\newpage}            %  Forced page break (new page).
\newtheorem{thm}{Theorem}
\newtheorem{lemma}[thm]{Lemma}

\newcommand{\cobord}{
$$\begin{picture}(140,130)  \footnotesize
    \put(-5,-5)       {\psfig{figure=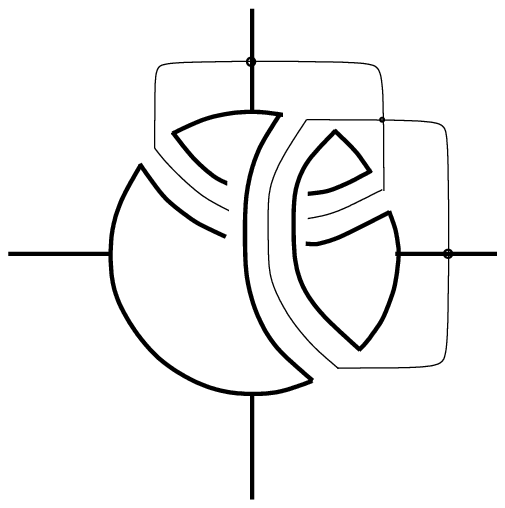}}
   \put(-10,57)       {$A$}
   \put(54,-3)         {$B$}
   \put(122,25)       {$A^t$}
   \put(25,123)       {$B^t$}
 \end{picture}$$}

\newcommand{\splitting}{
$$\begin{picture}(140,200)  \footnotesize
   \put(-225,-320)       {\psfig{figure=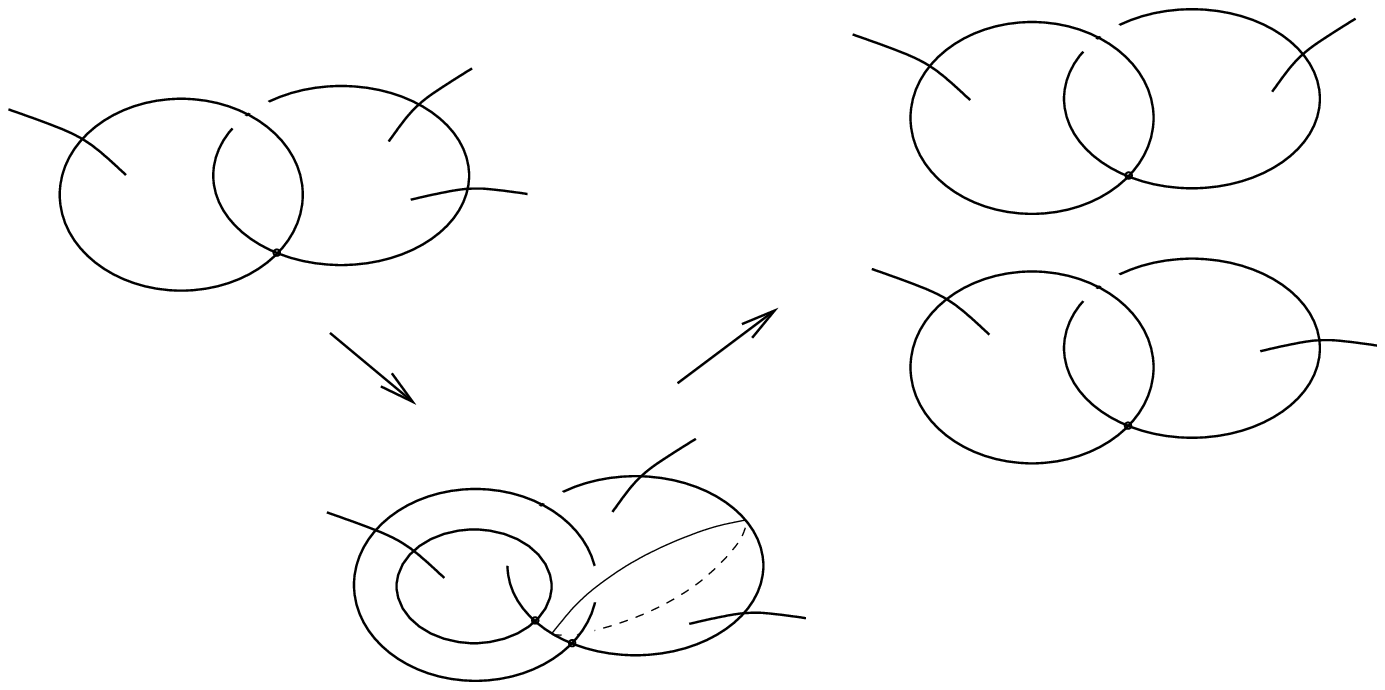}}
   \put(-107,168)       {$b$}
   \put(-95,108)       {$B$}
   \put(-10,180)       {$a'$}
   \put(13, 152)       {$a''$}
   \put(-25,112)       {$A$}
   \put(-20,52)       {$b''$}
   \put(-10,46)       {$b'$}
   \put(1, 22)       {$B'$}
   \put(1, -8)       {$B''$} 
   \put(55,68)       {$a'$}
   \put(88,12)       {$a''$}    \put(57,32)       {$\alpha$}
   \put(55,0)       {$A$}
   \put(140,190)       {$b'$}
   \put(130, 140)       {$B'$}
   \put(145, 122)       {$b''$}
   \put(128, 70)       {$B''$}
   \put(247, 192)       {$a'$}
   \put(247, 153)       {$A'$}
   \put(257, 108)       {$a''$}
   \put(243, 72)       {$A''$}
\end{picture}$$}

\newcommand{\whitney}{
$$\begin{picture}(140,160)  \footnotesize
   \put(-225,-330)       {\psfig{figure=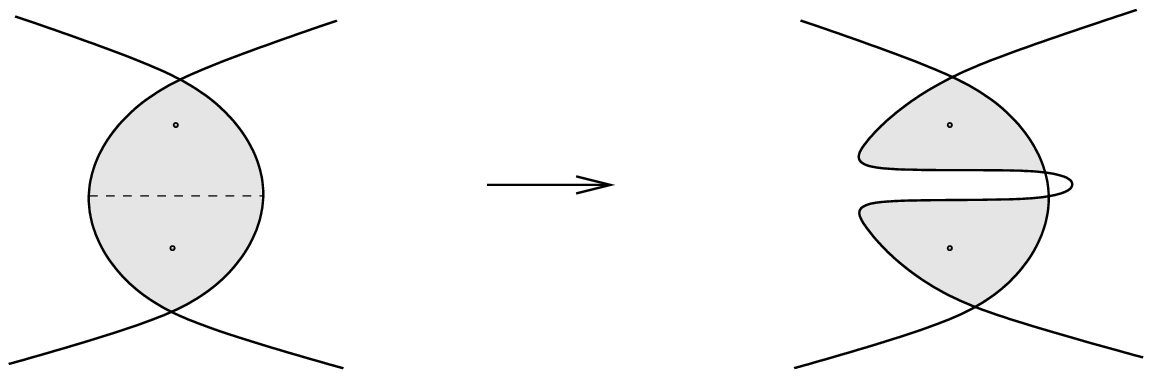}}
   \put(-95,30)       {$S_1$}
   \put(5,30)       {$S_2$}
   \put(-54,105)       {$a'$}
   \put(-55, 76)       {$a''$}
   \put(-32,80)       {$W$}
   \put(135,30)       {$S_1$}
   \put(235,32)       {$S_2$}
   \put(170, 107)       {$a'$}
   \put(170, 75)       {$a''$}  
   \put(187,102)       {$W'$}
   \put(188, 70)           {$W''$}
 \end{picture}$$}

\newcommand{\cycle}{
$$\begin{picture}(140,140)  \footnotesize
    \put(-240,-350)       {\psfig{figure=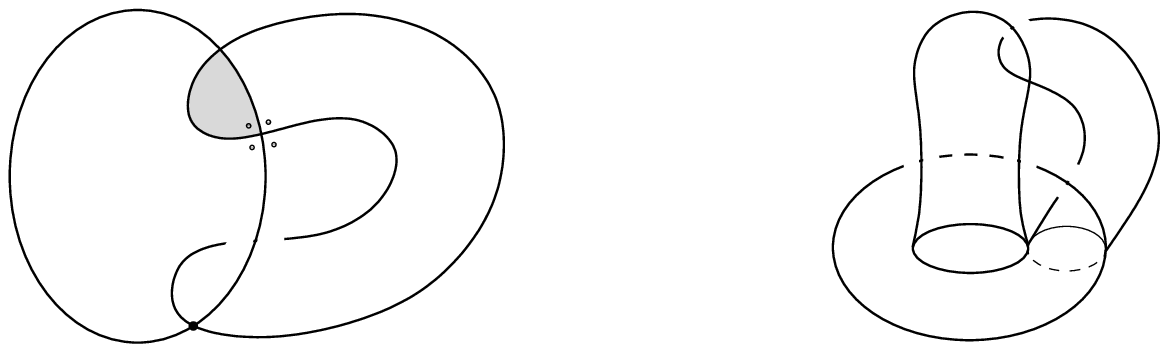}}
    \put(-94, 15)    {$B$}
    \put(-15, 15)    {$A$}
    \put(-52,92)   {$W$}
    \put(-28,72)    {$T$}
    \put(119, 8)     {capped torus $T^c$, dual to $W$}
 \end{picture}$$}

\newcommand{\unravel}{
$$\begin{picture}(140,140)  \footnotesize
    \put(-240,-410)       {\psfig{figure=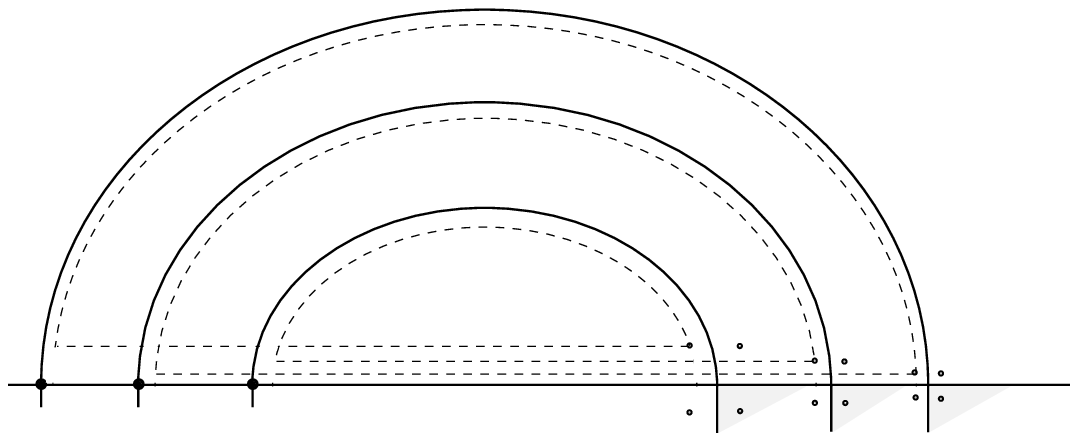}}
    \put(-57, 100) {$B_1$}
    \put(-37, 80) {$B_2$}    
    \put(-12, 60) {$B_3$}
    \put(115, 36) {$T_3$}
    \put(149, 31) {$T_2$}
    \put(178, 28) {$T_1$}
    \put(220, 15) {$A$}
    \put(124,2)     {$W_3$}
    \put(155,2)     {$W_2$}
    \put(182,2)     {$W_1$}
 \end{picture}$$}

\newcommand{\unravell}{
$$\begin{picture}(140,120)  \footnotesize
    \put(-235,-370)       {\psfig{figure=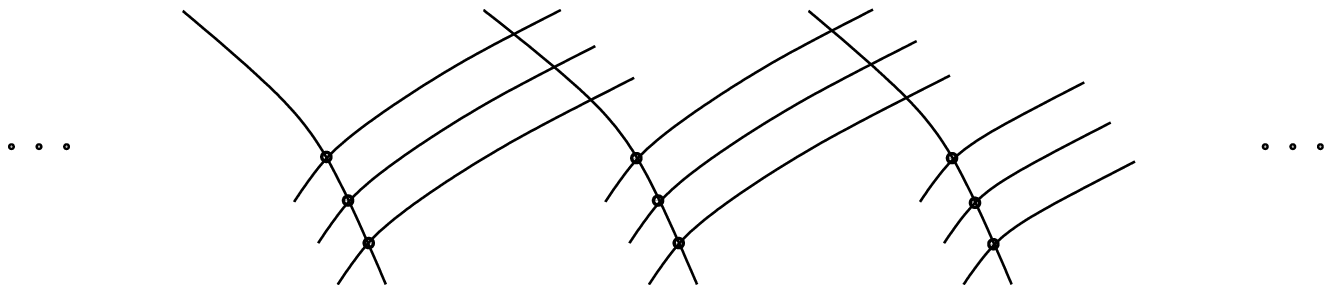}}
    \put(-89, 95) {$A^1$}
    \put(-60, 35) {$B^1_1$}    
    \put(-51, 17) {$B^1_2$}
    \put(-42, 2) {$B^1_3$}
    \put(-1,95) {$A^2$}
    \put(30, 35) {$B^2_1$}
    \put(39, 17) {$B^2_2$}
    \put(48,2){$B^2_3$}
    \put(95,95)     {$A^3$}
 \end{picture}$$}

\title[On disk embedding up to $S$-cobordism]{On disk embedding up to ${\mathbf S}$-cobordism}
\author{Vyacheslav S. Krushkal}
\thanks{Partially supported by NSF grant DMS 00-72722}
\address{Department of Mathematics, University of Virginia, Charlottesville, VA 22904}
\email{krushkal\char 64 virginia.edu}

\begin{document}

\begin{abstract}  
The disk embedding lemma is a technique underlying the topological
classification results in $4$-manifold topology for good fundamental groups. 
The purpose of this paper is to develop new tools for disk embedding that work
up to $s$-cobordism, without restrictions on the fundamental group. 
As an application we show that a surgery problem
gives rise to a collection of capped gropes that fit in the framework
of control theory.
\end{abstract}

\maketitle

The disk embedding lemma is a technique underlying the topological
classification results in $4$-manifold topology. 
It is conjectured to fail, in general, when the fundamental group 
of the $4$-manifold is free. 
A weaker lemma, disk embedding up to $s$-cobordism, would imply
$4$-dimensional surgery but not the $5$-dimensional $s$-cobordism 
conjecture. Recall its statement:

\vspace{.3cm}

\noindent
{\bf Conjecture}. 
{\sl Suppose $A$ is an immersion of a $2$-sphere into a $4$-manifold $M$, and there
is a framed immersed $2$-sphere $B$ with trivial algebraic self-intersection,
and algebraic intersection $1$ with $A$. Then there is an $s$-cobordism of
the immersion $A$ to a topological embedding.}

\vspace{.3cm}

The advantage of this formulation is the extra flexibility
in changing the ambient $4$-manifold, which could make it applicable in 
contexts where the disk embedding lemma in unavailable.
We refer the reader to \cite{Q1} for a discussion of the consequences of this
conjecture. It has been established \cite{FQ} under an additional assumption
that the transverse pair of spheres is ${\pi}_1$-null in $M$ (in other words,
the inclusion $A\cup B\longrightarrow M$ induces the trivial map on
fundamental groups.) 

The purpose of this paper is to develop new techniques for disk embedding that 
work up to $s$-cobordism, for any fundamental group. After reviewing the
foundational material -- the construction of $s$-cobordisms and splitting of gropes --
the following constructions are presented here. Sections \ref{splitting spheres}, 
\ref{splitting towers}: the grope splitting technique \cite{K}, \cite{KQ} is extended, 
up to $s$-cobordism, to transverse pairs of spheres, and more generally to pairs of 
Whitney towers. Section \ref{grope scobordism}: the $s$-cobordism construction for 
transverse pairs of spheres \cite{FQ} is adapted to the setting of capped gropes, so 
it can be applied at every surface stage. Section \ref{cycles}: A technique for 
eliminating cycles, which translates the algebraic-combinatorial data into a geometric 
configuration of gropes. In particular, the following statement is proved.

\vspace{.2cm}

\begin{thm} \label{geometric tree}
\sl Given a surgery problem, there is a procedure for converting it,
up to an $s$-cobordism, into a large collection of capped gropes whose
intersections (measured in ${\pi}_1 M^4$) are encoded by a tree, up 
to any given scale.
\end{thm} 

\vspace{.2cm}

This is a refinement of the outcome of grope splitting in \cite{KQ}.
In addition to the algebraic manipulations of the group elements
represented by the double point loops of the gropes (which were 
sufficient for the proof of the disk embedding lemma in the 
subexponential growth case in \cite{KQ}), this enables one to use 
the geometric techniques of control theory (cf \cite{FQ}, chapter 5.4.)

\vspace{.2cm}

{\em Acknowledgements.} I would like to thanks Frank Quinn for many discussions.

\section{Construction of ${\mathsc s}$-cobordisms} \label{construction}

The purpose of this section is to review the construction of $s$-cobordisms
in \cite{FQ}, Chapter 6. We include it here for convenience of the reader,
and also to fix the terminology, since this will be used throughout the paper.
Recall that this construction is used in \cite{FQ} to 
prove that a ${\pi}_1$-null immersion of a union of transverse pairs,
with algebraically trivial intersections, is $s$-cobordant to an embedding.

The starting point for the construction is a transverse pair of framed spheres
$(A, B)$ in a $4$-manifold $M$. Typically in the applications the transverse
pairs have algebraically trivial intersections. This means that $A$, $B$ 
have a ``preferred'' intersection point, and all other intersection and 
self-intersection points are paired up with Whitney disks.

\begin{figure}[ht]
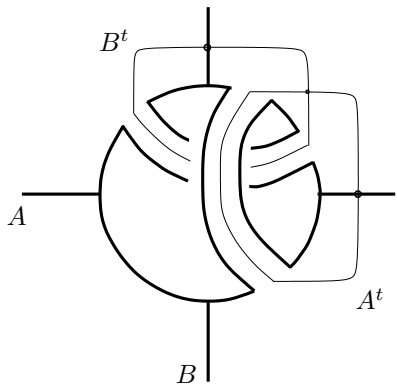

\small
\cobord
\caption{$s$-cobordism construction: transverse spheres}
\label{cobord}
\end{figure}

Let $D$ denote a ball in $M$ around the preferred intersection point
between $A$ and $B$, and consider the $4$-manifold $M'=M\smallsetminus D$. 
The intersection of $A$ and $B$ with $\partial D$ is the Hopf link $H$. 
Consider parallel copies $A'$, $B'$ of the spheres,
and let $S_1, S_2$ denote the circles of their intersection with $\partial D$.
$S_1$, $S_2$ are parallel copies of the components of $H$. Attach $5$-dimensional
$2$-handles to $M\times I$ along $S_1\times D^3$, $S_2\times D^3$ -- thus 
surgering $M$ along these circles. Figure~\ref{cobord} shows
the effect of surgery on $M'$. The cores of the handles (viewed in surgered $M'$)
are capped off with the disks bounded by $S_1$, $S_2$ in $\partial D$, 
slightly pushed out of the ball $D$.

This gives two spheres $A^t$, $B^t$, transverse to $A$, $B$ respectively.
They are embedded, but intersect in  
a single point, figure \ref{cobord}. The strategy for using them
in \cite{FQ} is the following. Consider the Clifford linking tori, one
for each pair of ``extra'' double points of $A$ and $B$. These tori
are geometrically dual to the Whitney disks for $A, B$, and they have caps 
provided by $A^t, B^t$. When $A\cup B$ is assumed to be ${\pi}_1$-null
in $M$, the resulting capped surfaces are used embedded
Whitney disks, thus solving the embedding problem for $A, B$.
The crucial observation is that the cores of the $2$-handles attached
to $M$ may be capped off by disks parallel to $A,B$, which are now embedded.
The resulting embedded $2$-spheres intersect the respective co-cores 
in a single point, and are used to attach $3$-handles.
The union of $M\times I$ and $2$- and $3$-handles is an $s$-cobordism
since the boundary map computing the relative homology with
${\mathbb{Z}}{\pi}_1$-coefficients is given by the identity matrix.

\vspace{.8cm}

\section{Grope Splitting} \label{grope splitting}

We refer the reader to \cite{FQ} for the definition and basic properties
of capped gropes. The (iterated) splitting operation from \cite{KQ} will be used
later in the paper, and is briefly reviewed in this section.

A disk-like grope is {\it dyadic\/} if all component surfaces are either disks 
or punctured tori. This means each non-cap surface has exactly one 
pair of dual subgropes attached to it. A grope has {\it dyadic 
branches\/} if all of the subgropes above the base level are dyadic. 

In a grope with dyadic branches, the caps are in one-to-one correspondence with
the {\em dyadic labels}: For each dual pair of surfaces in the grope, label
one by $0$ and the other by $1$. A cap gets a label by reading off the sequence of 
$0$s or $1$s encountered in a path going from the base surface to the cap.

The splitting operation splits a surface into two pieces, at the cost of 
doubling the dual subgrope. It can be used to decompose branches into 
dyadic branches, and can separate intersection points distinguished 
by properties unaffected by the dual doubling. In particular, it is used
to separate double points that have different dyadic labels, or 
that represent different group elements in ${\pi}_1 (M^4)$.

Suppose $A$ is a component surface of a grope, not part of the base. 
Let $B$ be the surface it is attached to, and $H$  the dual subgrope. 
Now suppose $\alpha$ is an embedded arc on $A$, with endpoints on the 
boundary, and disjoint from attaching circles of higher stages. 
In the 3-dimensional model, sum $B$ with itself by a tube about
$\alpha$ (the normal $S^{1}$ bundle), and discard the part of $A$ 
that lies  inside the tube. This splits $A$ into two components. 
$H$ is a dual for one component; obtain a dual for the other by taking a
parallel copy of $H$. The following lemma is proved by an inductive
application of this idea, descending from the caps downward to the first
stage surface.

\vspace{.2cm}

\begin{lemma} (Splitting) \label{just splitting} \sl 
Any grope can be transformed by iterated splitting to 
one with dyadic branches, and such that each cap satisfies:

\begin{enumerate}
\item there are no self-intersections;
\item all caps intersecting the given one have the same label; and 
\item the fundamental group classes of the loops through the 
intersection points are the same.
\end{enumerate}

Further, the subset of $\pi_{1}M$ occuring as double point loops is the same 
as that of the original grope.
\end{lemma} 

\vspace{.2cm}

We will need an iterated version of this lemma, where the intersections
are uniformized to distance $n$. This is made precise by the notion of an
$n$-{\em type}. We refer to \cite{KQ} for the formal definition of 
an $n$-type, and of a {\em collision} at distance $n$. 

\vspace{.2cm}

\begin{lemma} \label{n splitting} (Splitting to distance $n$) \sl 
Given a positive integer $n$ and a capped grope, there is a splitting so that 

\begin{enumerate}
\item every branch has an $n$-type,
\item there are no collisions at distance $\leq n$.
\end{enumerate}
\end{lemma}  

\vspace{.2cm}

Given an integer $n$ and a capped grope, its splitting to distance $n$
(outcome of lemma \ref{n splitting}) has the following algebraic-combinatorial
description. Consider any cap of the split grope. Its intersections
going out through chains of $n$ branches are encoded by a tree, whose valence
equals the class of the grope, and whose edges are labelled by a finite collection
of group elements represented by the intersections of the original grope.

The paths in this tree are mimicked by geometric moves on the caps: the
resulting double point loop represents the group element equal to the product
of the labels along the path in the tree. However, this combinatorial
graph is not, in general, mirrored geometrically as a tree of gropes, due 
to the existence of cycles. This is the subject of section \ref{cycles}.

\np

\section{Splitting up to ${\mathsc s}$-cobordism of transverse pairs} \label{splitting spheres}

This is a preliminary step
showing that there is an analogue of splitting of capped surfaces (lemma 
\ref{just splitting} above) 
for transverse pairs of spheres, up to $s$-cobordism. 
The $s$-cobordism is completed when eventually the immersed 
spheres are homotoped to embeddings.
An extension to transverse pairs with higher Whitney disk data is
discussed in section \ref{splitting towers}.

Let $A, B$ be a framed transverse pair in $M^4$, and suppose $A$ 
has two intersection points $a'$, $a''$ that need to be separated. 
Consider two parallel copies $B'$, $B''$ of $B$, and let $x', x''$ 
be the distinguished intersections between $B'$ and $A$, and 
between $B''$ and $A$, respectively. Let $\alpha$ be a simple 
closed curve in $A$, dividing it into two disks $D'$, $D''$, 
so that $x', a'\in D'$ and $x'', a''\in D''$. Attach a $2$-handle 
to $M^4\times I$ along $\alpha\times D^3$, thus surgering the $4$-manifold 
along the curve $\alpha$. Now both disks $D'$, $D''$ are capped
with copies of the core of the attached handle, giving immersed
spheres $A'$, $A''$. This is similar to the argument in section 
6.3 of \cite{FQ}, which converts a disk-embedding problem into
a sphere problem.

\begin{figure}[h]
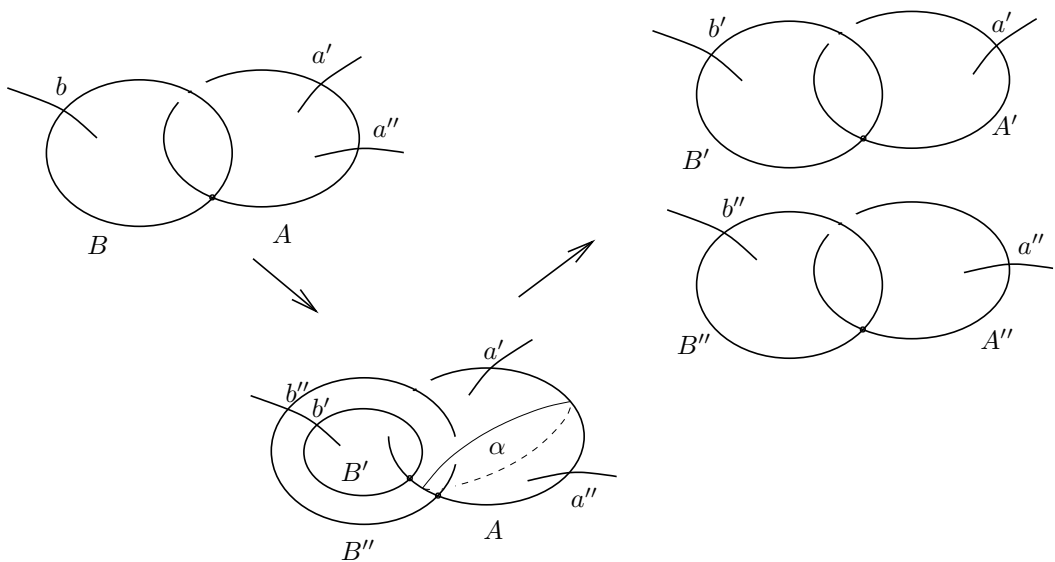

\small
\splitting
\caption{Splitting of a transverse pair}
\label{splitting}
\end{figure}

This creates two transverse pairs $(A', B')$ and $(A'', B'')$,
where $A'$ and $A''$ each inherits one of the intersections
$a'$, $a''$. The original intersection point $b$
is doubled, so both $B'$, $B''$ have an intersection.
If one of the spheres $A'$, $A''$ is homotopic to an embedding, 
then we will attach to it a $3$-handle, which together with 
the $2$-handle attached earlier yields an $s$-cobordism. 
To complete this $s$-cobordism, and also
to solve the original surgery problem, one needs to work with
both pairs of transverse spheres $(A', B')$, $(A'', B'')$.

As in lemma \ref{n splitting}, given any $n\geq 1$, there is an iterated splitting
of $(A,B)$ so that each sphere has an $n$-type, and there are no collisions
at distance $\leq n$. This is the analogue of splitting of capped surfaces.

\vspace{.8cm}

\section{Splitting of Whitney towers.} \label{splitting towers}
Here we show that there is an extension of the arguments above
from transverse pairs of spheres (dually: capped surfaces)
to transverse pairs together with layers of Whitney disks 
(dually: capped gropes). Let $(A, B)$ be a 
transverse pair with algebraically trivial intersections, 
so all extra intersection points are paired up by Whitney disks. Given $k\geq 1$, 
it may be arranged that: there are $k$ layers of Whitney disks for the
intersections of the spheres, and for each $1\leq i\leq k$ the interiors 
of Whitney disks at height $i$ are disjoint from all surfaces below 
that height (i.e. we have a transverse pair of Whitney towers
of height $k$.) This may be done, for example, by converting $(A,B)$ 
into a transverse pair of capped gropes of height $k$, and then 
contracting some of the surface stages. 

\begin{figure}[h]
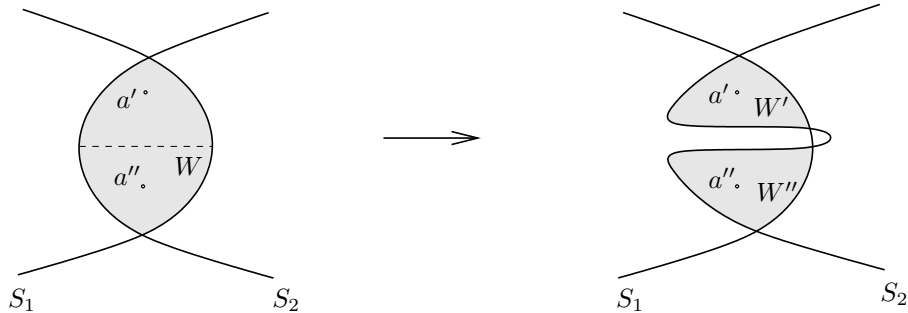

\small
\whitney
\caption{Splitting of a Whitney disk}
\label{splitting whitney}
\end{figure}

The splitting of transverse pairs described in section \ref{splitting spheres}
will be applied to {\em pairs} of intersections, rather than individual
intersections. Thus at all times the intersections will stay paired up,
and whenever a parallel copy of a surface (sphere or Whitney disk) is
taken, it acquires also a parallel copy of the whole corresponding 
higher Whitney data.

The remaining ingredient is the splitting of Whitney disks.
Suppose two intersections of some surfaces $S_1$, $S_2$ are paired up
with a Whitney disk $W$, and the interior of $W$ contains two
collections $a'$, $a''$ of points that need to be separated.
Consider a curve $\alpha$ in $W$ with one endpoint in $S_1$ and
the other one in $S_2$, separating $a'$ and $a''$. Perform a finger
move on one of the surfaces, say $S_1$, along $\alpha$, so that
two new $S_1$ -- $S_2$ intersections are created. The disk $W$
is replaced by two Whitney disks, $W'$ and $W''$, containing
$a'$, $a''$ respectively. This pushes the splitting problem to a lower
stage. Whitney towers are split, using an inductive application 
of this step together with the splitting of spheres described in section 
\ref{splitting spheres}.

\vspace{.8cm}

\section{Unraveling of short cycles.} \label{cycles}
Before introducing the general construction for capped gropes in section
\ref{grope scobordism}, we illustrate the problem, and our approach to 
it, in the setting of transverse pairs of spheres.
Start with a transverse pair of spheres with algebraically trivial
intersections, and split it to distance $n$. 
In this setup, the splitting of transverse pairs described above
will be applied to {\em pairs} of intersections, rather than individual
intersections. Thus at all times the intersections will stay paired up,
and whenever a parallel copy of a sphere is taken, its intersections 
acquire also a parallel copy of the corresponding Whitney disks.

\begin{figure}[h]
\small
\cycle
\caption{A cycle}
\label{cycle}
\end{figure}

Perform the $s$-cobordism construction (attach a pair of $2$-handles) 
in a neighborhood of the distinguished intersection point of each 
transverse pair of spheres, as described in section \ref{construction}. 
Now each Whitney disk has a dual capped torus.
We are going to focus on this collection $\{ T^c\}$ of capped tori, in the
complement of the original spheres. Algebraically we have a line, 
subdivided into intervals labelled by group elements, encoding the 
intersections of the spheres, as in \cite{KQ}. (Note that this line
also encodes, algebraically, the intersections among the $\{ T^c\}$.)

One cannot assume however that this picture holds 
geometrically as well, i.e. that there is a different capped torus 
for each vertex of the line. This is because of the presence of cycles.
By definition, a {\em cycle} is a path, embedded in the graph 
(in this context, in the line) whose endpoints correspond
to the same capped torus.
For example, suppose two spheres $A$, $B$ forming a transverse pair intersect
each other, as in figure \ref{cycle}, and let $g\in {\pi}_1 M$ be
the group element corresponding to the intersections. 

No splitting is necessary: both spheres have an $n$-type for each $n$, 
and the corresponding graph is a line, with each edge labelled by 
the same group element $g$. However, there is no line of tori
geometrically: the caps of the dual capped
torus intersect each other, this is a cycle of length $1$.

The following construction is used to unravel cycles, without
changing the uniform algebraic pattern of intersections. Before 
considering the general problem, here is the description for the
cycle in figure \ref{cycle}. Consider $n$ copies of one of the 
spheres, say $B$, denote them $B_1,\ldots, B_n$. Now there are
$n$ Whitney disks, $W_i$, $1\leq i\leq n$, and let $T_i$ denote
the corresponding Clifford torus, dual to $W_i$.

\begin{figure}[h]
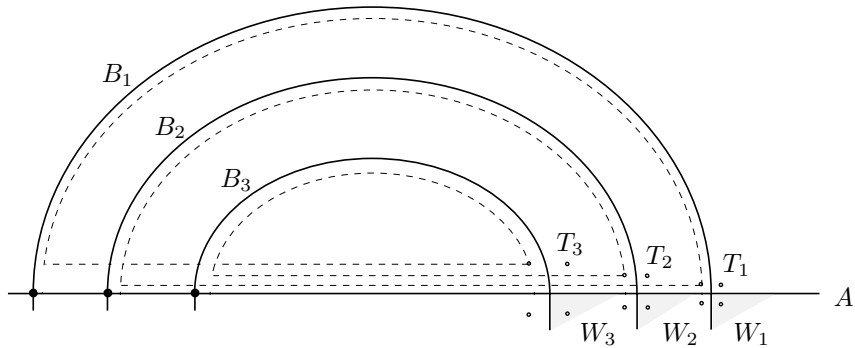

\small
\unravel
\caption{Unraveling of a cycle}
\label{unravel}
\end{figure}

Perform the construction of section \ref{construction} (attach a pair
of $2$-handles) in a neighborhood of the distinguished interesection 
point of each pair $(A, B_i)$, $1\leq i\leq n$. If all spheres
$A, B_1,\ldots, B_n$ are eventually shown to be homotopic to embedding,
their parallel copies will be used as attaching maps for $3$-handles. These,
together with the $2$-handles attached earlier, form an $s$-cobordism, 
since the boundary map $H_3\longrightarrow H_2$ is given by an upper 
triangular matrix.

Continue the construction of the capped tori: there is a unique choice
for the cap of $T^c_i$, intersecting $B_i$: it has to go over the
dual sphere to $B_i$, provided by the corresponding $2$-handle.
There are, however, $n$ choices for the other cap of each $T^c_i$,
intersecting $A$. We pick a cyclic shift, and it is this choice 
that is at the core of the
construction. For each $i=1,\ldots, n-1$ let the cap of $T^c_i$ go 
over the $2$-sphere, dual to $A$ near the intersection with $B_{i+1}$, 
and finally the cap of $T^c_n$ goes over the sphere dual to $A$
near the intersection with $B_1$. All caps are disjoint from each other, 
except for the unavoidable intersections that are inherited from the dual
spheres in the $s$-cobordism construction. Figure \ref{unravel} illustrates
the case $n=3$.

The outcome is a cyclic chain of capped surfaces $T^c_1,\ldots, T^c_n$,
with the caps of $T^c_i$ intersecting $T^c_{i-1}$ and $T^c_{i+1}$ --
the $n$-cyclic cover of the capped torus in figure \ref{cycle}.
This is a geometric realization, with injectivity radius $n/2$, of the 
algebraic graph encoding the intersections of the spheres.

We now continue the construction in the general case, for an arbitrary
transverse pair of spheres, split to distance $n$.
The solution here is going to be an implant of a segment of
length $n$ of the universal cover of figure \ref{unravel}.
More precisely, pick one of the transverse pairs, $(A^1,B^1)$, 
and look at the algebraic graph of intersections up to 
radius $n$. It corresponds to a sequence of pairs of spheres,
$(A^1,B^1), \ldots, (A^n, B^n)$ (actually, multiple parallel
copies due to splitting). We assign the labels $A_1, B_1, A_2, B_2,
\ldots$ to the spheres respecting their linear order. 
The edges of the graph correspond
to the intersections between $B^i$ and $A^{i+1}$. Moreover,
some of the pairs coincide (precisely in the presence of cycles.)

\begin{figure}[h]
\small
\unravell
\caption{}
\label{unravell}
\end{figure}

Recall that just like in the capped surface case, the spheres in
a pair are assigned dyadic labels. The notation may be misleading,
as $A$, $B$ do not correspond to such labels: $B^1$, $A^2$ may 
have the same dyadic label.

Replace $B^1$ with $n$ copies: $B^1_1,\ldots, B^1_n$. Since 
$B^1$ may occur elsewhere in the chain, say $B^1=B^i$ for some $i$,
all such $B^i$ are replaced with $n$ copies as well. 
Note that $B^i$ can never coincide with $A^j$ for any $i$, $j$, as this 
would imply
a collision at distance $<n$, and would contradict the splitting
assumption. If $B^2$ has not been affected, replace it with
$n$ copies. Continue this procedure until each $B^1,\ldots, B^n$
is replaced with $n$ copies. 

Now we have a segment of the infinite cyclic cover of figure
\ref{unravel}, and we follow the same construction as in that
example: assign caps
to the capped tori with a cyclic shift. The claim is that the result
is a sequence of length $n$ of capped tori, {\em without cycles},
exactly mimicing the algebraic
line of intersections of the spheres. If there were no cycles
to begin with, our construction just replicates the chain 
of capped surfaces, without altering the algebraic intersections. 
If there were cycles, they are eliminated.

\vspace{.8cm}

\section{${\mathsc s}$-cobordism construction for capped gropes} 
\label{grope scobordism}

This is the analogue of the $s$-cobordism construction for transverse pairs
of spheres (desrcibed in section \ref{construction}.) The new feature of the 
construction for capped gropes is that it can be applied at every surface stage.
At the end of this section we give the proof of theorem \ref{geometric tree}.

Recall that every surface (including the caps), above the first stage, in a capped
grope has a transverse grope. It is constructed using two copies of the dual
surface, cf section 1.4 in \cite{KQ}. To fix notations, suppose $S$ is a surface
in a capped grope, and let $A$, $B$ be surfaces attached to a symplectic
pair of circles ${\alpha}, {\beta}$ in $S$. Then the transverse capped grope $g_a$ for 
$A$ is built of two parallel copies of $B$ (and everything attached to $B$), 
and analogously the transverse grope $g_b$ for $B$ is made of two copies for $A$. The
transverse gropes $g_a, g_b$ are used to resolve intersections of surfaces
with $A$ and $B$ (and with other surface stages and caps attached to $A$, $B$) 
respectively.

The base surfaces $S_a$, $S_b$ of $g_a$, $g_b$ intersect in two points near the 
intersection point $p$ of the circles ${\alpha}, {\beta}$ in $S$. Consider 
a $4$-ball $D$ around the intersection point $p$ in $S$. The ball $D$ is chosen
so that $D_a=D\cap S_a$ and $D_b=D\cap S_b$ are disks, and $D$ contains the two 
intersection points $S_a\cap S_b$. Consider the circles $C_a$, $C_b$ of 
intersection of $S_a$, $S_b$ with $\partial D$. As in the construction of
section \ref{construction}, attach two $2$-handles to $M\times I$ along 
$C_a\times D^3$, $C_b\times D^3$. The cores of the $2$-handles (viewed in 
the surgered $4$-manifold) are capped with the disks $D_a$, $D_b$. The result 
is a pair of spheres, geometrically dual to the surfaces $S_a$, $S_b$. 
Each of these two spheres is embedded, but they intersect in two points. 

Consequently, the spheres may be used to resolve intersections with $A$, $B$
and with other surfaces attached to $A$, $B$, but the price is the presence
of new intersections between the surfaces pushed off $A$ and $B$.
These spheres are useful for finding embedded caps for the grope.
The crucial observation is that once the caps of the grope are improved to
embeddings, then the $3$-handles can be attached to complete the $s$-cobordism. 
The attaching $2$-spheres are constructed as follows: a parallel copy
of the core of each $2$-handle is capped off by an embedded disk built of two parallel
copies of the dual capped grope (an embedded disk may be found now,
since the caps are assumed to be embedded.) Attaching the $3$-handles completes
the construction of an $s$-cobordism, since the boundary map $C_3\longrightarrow C_2$
computing the relative homology with ${\mathbb{Z}}{\pi}_1$-coefficients is the
identity matrix.

We conclude with the proof of theorem \ref{geometric tree}. Starting with a 
transverse pair of spheres, convert one of them into a capped grope -- a capped
grope of height $2$ suffices, and we assume for simplicity that this is
the case. Let $n$ be a positive integer. Split the grope
to the distance $n$ (lemma \ref{n splitting}.) Fix a cap $C$ of the split grope;
as explained in section \ref{grope splitting}, the intersections among the caps 
of the split grope are encoded algebraically by a tree $T$, up to the distance $n$ 
from $C$. The goal is to convert this into a geometric picture where the capped
gropes correspond to the vertices of the tree, and the intersections
between their caps are encoded by the edges.
Consider all genus $1$ pieces of the base surface (branches) that appear 
in the tree $T$. 

The gropes have two surface stages and caps. The split grope has dyadic branches; 
fix a genus one piece of the base surface. It branches into two second 
stage surfaces, and then into four caps. Adapting to the setting of gropes 
the cyclic construction of section \ref{cycles}, consider $n$ parallel
copies of one of the second stage surfaces, and then $n$ parallel copies of
one in each dual pair of caps. The $s$-cobordism construction described above
will be applied at each intersection point of the attaching circles of all
surfaces (second stage surfaces and caps) present in the picture. This is the
model that will be applied at each branch of the split capped grope. 
Finally, the construction for resolving cycles in section \ref{cycles} 
is generalized from the the context of transverse pairs of spheres 
(intersections encoded by a line) to capped gropes (intersections
encoded by a tree.) The intersections between the caps are pushed down
and off the grope using the transverse spheres provided by the $s$-cobordism
construction. This is done with cyclic shift with respect to the numbering
of the parallel copies of the caps, and of the dual spheres.
The stage, to which the intersection is pushed down to, is determined
based on the dyadic labels of the caps. The proof that the cycles are resolved
is analogous to the sphere case, in section \ref{cycles}.


\vspace{.8cm}


\begin{thebibliography}{10}
\setlength{\parskip}{.1cm}

\bibitem[F]{F} M.H. Freedman, {\em Poincar\'{e} transversality and four-dimensional
surgery}, Topology 27 (1988), 171-175.

\bibitem[FQ]{FQ} M. Freedman and F. Quinn, {\em The topology of 
4-manifolds}, Princeton Math. Series 39, Princeton, NJ, 1990.

\bibitem[K]{K} V. Krushkal, {\em Exponential separation in $4$-manifolds},
Geom. Topol. 4 (2000), 397-405.

\bibitem[KQ]{KQ} V. Krushkal and F. Quinn, {\em Subexponential groups in 
$4$-manifold topology}, Geom. Topol. 4 (2000), 407-430.

\bibitem[Q]{Q} F. Quinn, {\em Ends of maps, III. Dimensions $4$ and $5$.}
J. Differential Geom. 17 (1982), 503-521.

\bibitem[Q1]{Q1} F. Quinn, {\em Problems in low-dimensional topology}, Annals Math.
Studies 149, 423-436. Princeton University Press, 2001.


\end{thebibliography}
\end{document}